
\magnification=900
\baselineskip=13.4pt

\overfullrule=0pt

\def\EEE{{\bf E}}       
      
      \def\PPP{{\bf P}} 
  \def\RRR{{\bf R}}     
      
  \def\ZZZ{{\bf Z}}

\def\AA{{\cal A}}       
\def\EE{{\cal E}}  \def\FF{{\cal F}}     
      \def\LL{{\cal L}}

\def\today{\ifcase\month\or 
  January\or February\or March\or April\or  
  May\or June\or July\or August\or  
  September\or October\or November\or
  December\fi\space\number\day,\ \number\year}

\def\ref#1{{\rm [}{\bf #1}{\rm ]}}   
\def\nref#1#2{{\rm [}{\bf #1}{\rm ;\ #2]}}

\def\card{{\mathop{\rm card}}}

\def\b#1{\nin{\bf#1.}}

\outer\def\proclaim#1{\medbreak\noindent\bf\ignorespaces
   #1\unskip.\enspace\sl\ignorespaces}
\outer\def\endproclaim{\par\ifdim\lastskip<\medskipamount\removelastskip
   \penalty 55 \fi\medskip\rm}

\def\nin{\noindent}

\def\comp{{\leavevmode
     \raise.2ex\hbox{${\scriptstyle\mathchar"020E}$}}}

\font\tenBbb=msbm10

\newfam\Bbbfam \textfont\Bbbfam=\tenBbb
\font\sevenBbb=msbm7
\scriptfont8=\sevenBbb

\def\rect#1#2#3{\raise .1ex\vbox{\hrule height.#3pt
   \hbox{\vrule width.#3pt height#2pt \kern#1pt\vrule width.#3pt}
        \hrule height.#3pt}}

\def\qed{$\hskip 5pt\rect364$} 

\def\1111{\bf 1}

%
%

\centerline{\bf   Kemeny's Constant for Markov Processes}
\medskip
\centerline{P.J. Fitzsimmons}
\centerline{Department of Mathematics, UC San Diego}
\centerline{La Jolla, CA 92093-0112}
\centerline{USA}
\centerline{\tt pfitzsim@ucsd.edu}
\medskip
\centerline{\today}

\footnote{}{AMS 2020 Mathematics Subject Classication: Primary 60J10; Secondary 60J25, 60J45, 60J55.}
\footnote{}{Keywords and phrases: Kemeny's constant; hitting time; time reversal; mean occupation measure; duality; Markov chain;  Hunt process; positive recurrence; local time; potential theory; effective resistance metric.}
\medskip
\bigskip

\centerline{\bf ABSTRACT}
\medskip

The mean time taken by  an irreducible Markov chain on a finite state space to hit a target chosen at random according to the stationary distribution does not depend on the initial state of the chain. This mean time is known as Kemeny's constant.
I present a new approach, based on time reversal and a mean occupation time formula. 
 The method is used to prove an analogous result  for continuous-time  Markov processes. We also present a second approach, based on work of N.~Eisenbaum and H.~Kaspi, when all states are regular.
 Examples are provided.
 \bigskip \bigskip

\moveright2in\vbox{\hrule width 2in}

 \bigskip \bigskip

\def\EEE{{\bf E}}       
      
      \def\PPP{{\bf P}} 
  \def\RRR{{\bf R}}     
      
  \def\ZZZ{{\bf Z}}

\def\AA{{\cal A}}       
\def\EE{{\cal E}}  \def\FF{{\cal F}}     
      \def\LL{{\cal L}}

\def\today{\ifcase\month\or 
  January\or February\or March\or April\or  
  May\or June\or July\or August\or  
  September\or October\or November\or
  December\fi\space\number\day,\ \number\year}

\def\reef#1{{\rm [}{\bf #1}{\rm ]}}   

\def\nref#1#2{{\rm [}{\bf #1}{\rm ;\ #2]}}

\def\card{{\mathop{\rm card}}}

\def\b#1{\nin{\bf#1.}}

\outer\def\proclaim#1{\medbreak\noindent\bf\ignorespaces
   #1\unskip.\enspace\sl\ignorespaces}
\outer\def\endproclaim{\par\ifdim\lastskip<\medskipamount\removelastskip
   \penalty 55 \fi\medskip\rm}

\def\nin{\noindent}

\def\comp{{\leavevmode
     \raise.2ex\hbox{${\scriptstyle\mathchar"020E}$}}}

\font\tenBbb=msbm10

\newfam\Bbbfam \textfont\Bbbfam=\tenBbb
\font\sevenBbb=msbm7
\scriptfont8=\sevenBbb

\def\rect#1#2#3{\raise .1ex\vbox{\hrule height.#3pt
   \hbox{\vrule width.#3pt height#2pt \kern#1pt\vrule width.#3pt}
        \hrule height.#3pt}}

\def\qed{$\hskip 5pt\rect364$} 

\def\1111{\bf 1}

\def\AF{1} 
\def\AH{2} 
\def\AKDR{3} 
\def\BHLMT{4} 
\def\BG{5} 
\def\CFS{6} 
\def\C{7} 
\def\CG{8} 
\def\CZ{9} 
\def\Dev{10} 
\def\DMT{11}
\def\D{12} 
\def\EK{13} 
\def\F{14} 
\def\FG{15} 
\def\FK{16} 
\def\Gexc{17} 
\def\Gtrans{18} 
\def\Gexcess{19} 
\def\GK{20} 
\def\GS{21} 
\def\H{22} 
\def\HS{23} 
\def\HY{24} 
\def\HunMix{25} 
\def\Hun{26} 
\def\Iso{27}
\def\IM{28} 
\def\KM{29} 
\def\KS{30} 
\def\K{31} 
\def\LeLo{32} 
\def\LP{33}
\def\LiLy{34} 
\def\M{35} 
\def\MT{36} 
\def\Pin{37} 
\def\Sh{38} 
\def\Si{39} 
\def\W{40} 



\b{1. Introduction}
\medskip

Let $X$ be a positive recurrent Markov process with  state space $E$ and unique stationary distribution $\pi$. Suppose that for each pair $(x,y)$ of states we have $\PPP^x[T_y<\infty]=1$, where $T_y$ denotes the hitting time of $y$. It was shown in \reef{\KS}, in the context of Markov chains with finite $E$, that the ``mean time to equilibrium''
$$
K(x):=\int \EEE^x[T_y]\,\pi(dy),\qquad x\in E,
\leqno(1.1)
$$
does not, in fact, depend on the starting state $x$. Recently, Pinsky \reef{\Pin} has demonstrated the same result in the context of $1$-dimensional diffusions. For further background discussion, see pp.~1311--1312 in \reef{\Hun}, and also \reef{\Dev}.

In this paper we present a new  proof of the constancy of the ``Kemeny function'' $K$ defined above.
We will first prove the result in the case of discrete-time Markov chains as a way to highlight the main idea, shorn of technical issues. The proof may not be  as brief as Doyle's \reef{\D}, but the method can be made to work in a much more general  context. 

The Markov chain case  is the subject of section 2. Section 3 treats the case of continuous-time strong Markov processes in duality. The preliminary Theorem (3.9) establishes the constancy of $K$ a.s.\ with respect to the invariant distribution $\pi$.
 In section 4   we refine the basic result, showing that the exceptional set is in fact empty. This hinges on the apparently new result that $\sup_x\EEE^x[T_y]<\infty$ for all $y\in E$, in complete generality. Section 5 concerns the situation in which all states are regular. We use a different method, based on a symmetry noticed by Eisenbaum and Kaspi \reef{\EK}, to show that $K$ is constant.
 Section 6 contains four examples.
\medskip

\b{Notation} If $(F,\FF)$ is a measurable space we use $b\FF$ to denote the vector space of bounded $\FF$-measurable functions mapping $F$ into $\RRR$; likewise  $p\FF$ is the convex cone of  $\FF$-measurable functions mapping $F$ into $[0,\infty[$. And $bp\EE$ is just $b\EE\cap p\EE$. If  $\mu$ is a  measure on $(F,\FF)$ we sometimes write  $\mu(f)$ for $\int_F f\,d\mu$.
\bigskip

\b{2. Markov Chain}
\medskip

In this section   $X=(X_n)_{n\ge 0}$ is  a  discrete-time Markov chain with finite  state space $E=\{1,2,\ldots,N\}$ and one-step transition matrix $P$. We assume that $X$ is irreducible and use  $\pi$ to denote the unique stationary distribution for $X$. Thus $\pi P=\pi$ and $\pi\cdot{ \1111}=1$.  (Here ${\bf 1}$ is an $N\times 1$ column of $1$s.)

The law of $X$ started at $x\in E$ is $\PPP^x$, on the sample space $\Omega=E^{\{0,1,2,\ldots\}}$ of all $E$-valued sequences $\omega=(\omega_n)_{n\ge 0}$. We realize $X$ as   the coordinate process:  $X_n(\omega) = \omega_n$.   The  symbol $\EEE^x$ will  be used for  expectation based on  $\PPP^x$, and if $\mu$ is a probability measure on $E$ then $\PPP^\mu:=\sum_{x\in E}\mu(x)\PPP^x$ denotes the law  of $X$ under the initial distribution $\mu$; $\EEE^\mu$ is the associated expectation.

\def\hP{\hat P}

Because $\pi$ is invariant for $P$, the recipe
$$
\hP_{xy} := P_{yx}{\pi_y\over\pi_x},\qquad x,y\in E,
\leqno(2.1)
$$
defines a stochastic matrix $\hat P$  (the $\pi$-dual of $P$) that is also irreducible and that admits $\pi$ as unique stationary distribution.
The Markov chain $\hat X=(\hat X_n)_{n\ge 0}$ associated with $\hP$ is realized as the coordinate process on $\Omega$, governed by laws 
 $\hat\PPP^x$, $x\in E$. We write $\hat X_n(\omega) =\omega_n$ when referring to $\hat X$.
  By (2.1),
$$
\pi_x\,\PPP^x[X_1=y]=\pi_y\,\hat\PPP^y[\hat X_1=x],
$$
and then by repeated application of the Markov property:
$$
\eqalign{
  \pi_x\,\PPP^x[X_1=x_1,& X_2=x_2,\ldots , \,X_{n-1}=x_{n-1}, X_n=y]\cr
&=\pi_y\,\hat\PPP^y[\hat X_1= x_{n-1},\ldots,\hat X_{n-2} =x_2, \hat X_{n-1} = x_1, \hat X_n=x],\cr
}
\leqno(2.2)
$$
for all $x,y,x_x,\ldots,x_{n-1}$ in $E$ and $n\ge 1$. In this sense, $\hat X$ is $X$ with the direction of time reversed.

Let $D_z:=\min\{n\ge 0: X_n=z\}$ denote the entry time into the state $z$. It is well known that the tail $n\mapsto \PPP^x[D_z >n]$ decays geometrically; in particular the expectations $\EEE^x[D_z]$ are all finite. Define the ``Kemeny function'' $K$, for $x\in E$, by
$$
K(x):=\sum_{z\in E}\EEE^x[D_z]\pi_z=\EEE^x[D_Z]<\infty,
\leqno(2.3)
$$
where $D_Z$ denotes the ``equilibrium time''; that is, the entry time to a state  $Z$ chosen randomly (and independently of $X$) using $\pi$. 
[To be more rigorous about the second equality in (2.3), one may define $D_Z$ on $\overline\Omega:=\Omega\times E$ as $D(\omega,z):=D_z(\omega)$, endowing $\overline\Omega$ with the laws $\PPP^x\otimes \pi$, $x\in E$. The dual  time $\hat D_Z$ is treated the same way; details are left to the reader.]
\medskip

\proclaim{(2.4) Proposition \nref{\KS}{Thm.~4.4.10}} The function $K$ is constant on $E$. In fact,
$$
K(x) =\hat\EEE^\pi[\hat D_Z],\qquad\forall x\in E.
\leqno(2.5)
$$
\endproclaim

By the evident symmetry of our hypotheses,  the dual equality $\hat K(y)=\EEE^\pi[D_Z]$ holds as well. In particular, the two Kemeny constants 
$$
\sum_{x,y\in E} \pi_x\pi_y\,\EEE^x[D_y],\qquad\sum_{x,y\in E}\pi_x\pi_y\,\hat\EEE^y[\hat D_x]
$$
are equal, and the functions $K$ and $\hat K$ coincide.

Our proof  of Proposition (2.4)  is based on the time-reversal duality (2.2) and  a  mean occupation time formula \nref{\AF}{Prop.~2.4} due  to Aldous and Fill,  Lemma (2.6) below.  The idea behind such a result goes back to Chung \nref{\C}{Thm.~2;  Appendix} in discrete time, and was used by Harris \nref{\H}{Thm.$\,$1}  (discrete time) and Silverstein \nref{\Si}{Thm.$\,$9} (continuous time) as a way to construct invariant measures.  In the statement of the lemma we use the notation  $N_n(y):=\sum_{k=0}^{n-1} 1_{\{X_k=y\}}$ for  the number of visits to $y$ before time $n$. 

\proclaim{(2.6) Lemma} Let $\mu$ be a probability distribution on $E$ and let $S$ be a (possibly randomized) stopping time of $X$ such that $\EEE^\mu[S]<\infty$ and $\PPP^\mu[X_S = y] =\mu(y)$, for all $y\in E$. Then
$$
\EEE^\mu[N_S(y)]=\pi_y\,\EEE^\mu[S],\qquad\forall y\in E.
\leqno(2.7)
$$
\endproclaim

\nin We will state and prove a continuous-time version of this result in section 3.

\medskip

\nin{\sl Proof of Proposition (2.4).} 
Because $\{n<D_z\} =\cap_{k=0}^n\{X_k\not=z\}$,   (2.2)  implies  that
$$
\pi_x\,\PPP^x[X_n=y, n<D_z]=\pi_y\,\hat\PPP^y[\hat X_n=x,n<\hat D_z],\qquad x,y\not=z, n=1,2,\ldots.
$$
Summing on $n\ge 0$:
$$
\pi_x\,\EEE^x[N_{D_z}(y)]=\pi_y\,\hat\EEE^y[\hat N_{\hat D_z}(x)],\qquad \forall x,y\in E.
\leqno(2.8)
$$
(Both sides vanish if either $x=z$ or $y=z$.) Now sum over $y\in E$ to get
$$
\pi_x\,\EEE^x[D_z] =\hat\EEE^\pi\,[\hat N_{\hat D_z}(x)],\qquad \forall x\in E.
\leqno(2.9)
$$
Finally, multiply both sides of (2.9) by $\pi_z$ and sum   to obtain
$$
\pi_x\,\EEE^x\,[D_Z]=\hat\EEE^\pi[\hat N_{\hat D_Z}(x)].
\leqno(2.10)
$$
We now apply     Lemma (2.6)  (for $\hat X$) to rewrite the right side of (2.10):
$$
\hat\EEE^\pi[\hat N_{\hat D_Z}(x)]=\pi_x\,\hat\EEE^\pi[\hat D_Z],\qquad\forall x\in E.
\leqno(2.11)
$$
Taken together, (2.10) and (2.11) prove   Proposition (2.4).\qed
\bigskip

\b{(2.12) Remark} It is worth noting that the first return times $T_z:=\min\{n\ge 1: X_n=z\}$ enjoy an analogous identity. To see this note that $\EEE^x[T_z] =\EEE^x[D_z]$ if $x\not=z$, while
$$
\EEE^z[T_z] = 1+\sum_{y\in E} P_{zy}\EEE^y[D_z],
$$
so that, letting $T_Z$ denote the ``$\pi$ mixture'' of the $T_z$, we have
$$
\EEE^x[T_Z] = 1+\EEE^x[D_Z],
\leqno(2.13)
$$
because $\pi_x=1/\EEE^x[T_x]$.
\bigskip

\b{3. Continuous-time  Markov Process}
\medskip

For the rest of the paper we take $X$ to be a continuous-time Hunt process  $X=(X_t)_{t\ge 0}$ with  Lusin state space $(E,\EE)$. That is, $X$ is a Borel right Markov process with quasi-left continuous sample  paths. Good  references for such processes are \reef{\BG} and \reef{\Sh}. We take $X$ to be the coordinate process defined on the sample space $\Omega$ of right-continuous, left-limited paths from $[0,\infty[$ to $E$. As before, $\PPP^x$ is the law on $\Omega$ for $X$ started in state $x\in E$; and the notations $\EEE^x$, $\EEE^\mu$, etc., are as before.

Let $(P_t)_{t\ge 0}$ denote the transition semigroup for $X$. We assume that $X$ is {\it honest} in the sense that $P_t 1(x) =1$ for all $x\in E$ and all $t\ge 0$.
To be able to define $K$ in the present context we assume that $X$ hits points:
$$
\PPP^x[T_y<\infty]=1,\qquad\forall x,y\in E,
\leqno(3.1)
$$
where $T_y:=\inf\{t>0: X_t=y\}$ denotes the hitting time of $y$. It follows from (3.1) that  the only polar set for $X$  is the empty set. Moreover, by \nref{\KM}{Thm.$\,$1} (take $\nu$ there to be the point mass at some fixed point of $E$) the process $X$ is Harris recurrent. Thus $X$ admits a unique (up to a multiple)  invariant measure $\xi$; see \nref{\AKDR}{Thm.$\,$I.3}. We {\it assume} $\xi(E)<\infty$ and then normalize $\xi$ to obtain the  stationary distribution $\pi$ for $X$, so that  $\pi P_t = \pi$ for all $t\ge 0$ and $\pi(E)=1$.
 Thus,  $X$ is {\it positive recurrent\/}. More precisely, $\pi$ is conservative \nref{\Gexcess}{p.$\,$8; Cor.$\,$(3.9)}, and unique in the sense that  any other excessive measure for $X$ is proportional to $\pi$; for this use \nref{\HY}{Thm.$\,$2.1}, noting that (3.1) implies that $X$ is ``finely irreducible.''

In particular,  $\pi$ is a reference measure: for all  $A\in\EE$,
$$
\pi(A) = 0\qquad\Longrightarrow\qquad \EEE^x\left[\int_0^\infty 1_A(X_t)\, dt\right]=0,\forall x\in E.
\leqno(3.2)
$$
[To see this pick a state $y$. If $y$ is irregular then
$$
\EEE^y\int_0^{T_y} 1_A(X_t)\,dt
$$
defines a $\sigma$-finite measure ({\it cf.$\!$} (4.1) below) that is invariant (hence proportional to $\pi$) and also charges all sets of positive potential. If $y$ is regular, the same is true of
$$
A\mapsto c\cdot 1_A(y)+n_y\int_0^{T_y} 1_A(X_t)\,dt,
$$
where $c\ge 0$ is a certain constant and  $n_y$ is the It\^{o} excursion law for the excursions of $X$ from $y$. This is the construction alluded to in section 2, just above (2.6); see \nref{\Gexc}{Thm.$\,$(8.1)}.]

As a replacement for the explicit recipe (2.1) we {\it assume}  that $X$ admits a dual process $\hat X$ with respect to $\pi$. This is a  second Hunt  process $\hat X=(\hat X_t)_{t\ge 0}$ with the same state space  as $X$, and   such that the  transition semigroup $(\hat P_t)_{t\ge 0}$ of   $\hat X$ is dual to that of $X$:
$$
\int_E f(x)P_tg(x)\,\pi(dx) =\int_E \hat P_tf(y) g(y)\,\pi(dy),\qquad\forall t>0, \quad \forall f,g \in b\EE.
\leqno(3.3)
$$
We assume that $\hat X$ is honest and satisfies the analog of (3.1). In particular,  we are in the context of Chapter 6 of \reef{\BG}.
Clearly  $\pi$ is the unique invariant distribution for $\hat X$.

Because $\pi$ is a reference measure,  the resolvent operators $U^\alpha:=\int_0^\infty e^{-\alpha t} P_t\,dt$ and $\hat U^\alpha:=\int_0^\infty e^{-\alpha t} \hat P_t\,dt$, $\alpha>0$,  are absolutely continuous with respect to $\pi$; there is   a dual density $(x,y)\mapsto u^\alpha(x,y)$ such that
$$
U^\alpha f(x) =\int_E u^\alpha (x,y) f(y)\,\pi(dy),\qquad x\in E,
\leqno(3.4)
$$
and 
$$
\hat U^\alpha f(y) =\int_E u^\alpha (x,y) f(x)\,\pi(dx),\qquad y\in E,
\leqno(3.5)
$$
for all bounded or positive $\EE$-measurable $f$. Moreover, $x\mapsto u^\alpha(x,y)$ is $\alpha$-excessive for each $\alpha>0$ and $y\in E$; likewise $y\mapsto u^\alpha(x,y)$ is $\alpha$-co-excessive (that is, $\alpha$-excessive with respect to $\hat X$) for  each $\alpha>0$ and $x\in E$. See \nref{\BG}{Thm.$\,$VI(1.4)}.

 It will be useful later to know that because of (3.2),  $\pi$ charges each non-empty finely open (or co-finely open) set.

Our goal is to show that
$$
E\ni x\mapsto \int_E \EEE^x[T_z]\,\pi(dz), 
$$
is a constant function. In the present generality an issue arises that is not present in the context of section 2. For example, Hunter's estimate \nref{\HunMix}{Thm.~4.2} 
$$
\EEE^x[D_Z]\ge{\card(E)-1\over 2},
\leqno(3.6)
$$
for discrete-time chains with finite state space $E$,  implies fairly easily  that Kemeny's constant for a discrete-time Markov  chain with {\it infinite}  state space must be infinite. See \nref{\LiLy}{Thm.~2.1}, and   \reef{\AH}; see also \nref{\LeLo}{p.$\,$744} for (3.6) in the reversible case. Pinsky \reef{\Pin} has shown  that the mean time to equilibrium for a $1$-dimensional diffusion is a constant function of the starting point, but that the  constant is finite if and only if any boundary points are entrance boundaries. The parallel result for continuous-time birth-and-death processes on $E=\{0,1,2,\ldots\}$ can be found in \reef{\BHLMT}.

The following is probably  well known, and will be important in sections 4 and  5. The proof is postponed until section 4 where a more detailed discussion of these issues appears, including a uniform version of (3.8); see (4.10).

\proclaim{(3.7) Proposition} For all $x\in E$ and $y\in E$,
$$
\EEE^x[T_y]<\infty.
\leqno(3.8)
$$
\endproclaim

The finiteness of  $K(x)$ is no longer guaranteed (as it is when $E$ is finite). But see Lemma (3.16)(b) below for the relevant solidarity result.
We  write $T_Z$ for the randomized hitting time obtained by choosing a target state $Z$ from $E$ using $\pi$, independently of $X$. The dual time $\hat T_Z$ is defined analogously.

\proclaim{(3.9) Theorem} Define $K(x):=\EEE^x[T_Z]$ and $\hat\kappa:=\hat \EEE^\pi[\hat T_Z]$.  Then
$$
K(x) \le\hat\kappa,\qquad\forall x\in E,
\leqno(3.10)
$$
and
$$
K(x)=\hat\kappa,\qquad\hbox{for }\pi\hbox{-a.e. } x\in E.
\leqno(3.11)
$$
\endproclaim

The constant $\kappa$ is defined analogously as $\int_E K(x)\,\pi(dx)=\EEE^\pi[T_Z]$. Integrating in (3.11) we find that $\kappa = \hat\kappa$.

Our  proof of Theorem (3.9) depends on a continuous-time version of  Lemma (2.6). 

\proclaim{(3.12) Lemma} Let $\mu$ be a probability distribution on $E$ and let $S$ be a (possibly randomized) stopping time of $X$ such that $\PPP^\mu[S<\infty]=1$ and $\PPP^\mu[X_S \in\cdot\, ] =\mu$. Suppose that the measure (in $f$) defined by the left side of (3.13) below is $\sigma$-finite.  Then $\EEE^\mu[S]<\infty$ and 
$$
\EEE^\mu\left[\int_0^S f(X_t)\,dt\right]=\pi(f)\,\EEE^\mu[S],\qquad\forall f\in bp\EE.
\leqno(3.13)
$$
\endproclaim

\medskip

\nin{\sl Proof.} We  view the left side of (3.13) as $\eta(f)$ for a certain measure $\eta$. Clearly $\eta(1) =\EEE^\mu[S]$. We show that $\eta$ is invariant: $\eta P_t=\eta$ for all $t>0$. Fix $t>0$, and observe that $\eta P_t(f) =\eta(P_t f)$ is equal to
$$
\eqalign{
\EEE^\mu\int_0^S P_t&f(X_u)\,du \cr
&=\int_0^\infty \EEE^\mu\left[ P_tf(X_u); u<S\right]\,du\cr
&=\int_0^\infty \EEE^\mu\left[ f(X_{t}\comp\theta_u); u<S\right]\,du 
=  \EEE^\mu\left[ \int_t^{S+t}f(X_{v})\,dv \right]\cr
&= \EEE^\mu\left[ \int_0^{S}f(X_{v})\,dv \right]- \EEE^\mu\left[ \int_0^tf(X_{v})\,dv \right]+ \EEE^\mu\left[ \int_S^{S+t}f(X_{v})\,dv \right].\cr
}\leqno(3.14)
$$ 
(Here $\theta_u$ is the usual shift operator on the sample space of $X$.) By the strong Markov property the third term on the far right of (3.14) is equal to
$$
\EEE^\mu\left[ \left(\int_0^{t}f(X_{v})\,dv\right)\comp\,\theta_S \right]=\EEE^\mu\left[ \EEE^{X_S}\left[\int_0^{t}f(X_{v})\,dv\right]\right],
\leqno(3.15)
$$
which is equal to the (finite) term being subtracted in (3.14) because the $\PPP^\mu$ law of $X_S$ is $\mu$. As $t>0$ was arbitrary, this demonstrates the asserted invariance. In particular, being $\sigma$-finite, $\eta$ is an excessive measure for $X$. Consequently, $\eta=c\cdot\pi$ for some constant $c\in[0,\infty[$. Clearly $\EEE^\mu[S] =\eta(1)=c\cdot\pi(1) =c<\infty$. 
\qed
\medskip

We also need the following result concerning  mean hitting times and the Kemeny function.

\proclaim{(3.16) Lemma} (a) We have
$$
\EEE^x[T_z]\le\EEE^x[T_y]+\EEE^y[T_z],\qquad \forall x,y,z\in E.
\leqno(3.17)
$$
In particular,
$$
K(x)\le \EEE^x[T_y]+K(y),\qquad\forall x,y\in E.
\leqno(3.18)
$$

(b) Either $K(x)<\infty$ for all $x\in E$ or $K(x) =\infty$ for all $x\in E$.
\endproclaim
\medskip

\nin{\sl Proof.} (a)  Abbreviate $\tau=T_{y}$. Then by the {\it terminal time}  property of hitting times
$$
T_z = t+T_z\comp\,\theta_t \quad \hbox{ on }\quad \{t<T_z\},\qquad  \forall t>0,
$$
we have
$$
T_z=T_z\cdot 1_{\{T_z\le\tau\}}+(\tau+T_z\comp \theta_\tau)1_{\{T_z>\tau\}} = (T_z\wedge\tau) +T_z\comp \theta_\tau1_{\{T_z>\tau\}}.
$$
Consequently,
$$
\EEE^x[T_z]\le \EEE^x[\tau]+\EEE^x[\EEE^{X_\tau}[T_z]].
$$
which implies (3.17) because $\PPP^x[X_\tau = y,\tau<\infty] =1$.
Integrating out $z$ with respect to $\pi$ yields (3.18).  Part (b) follows immediately from (3.18) and Proposition (3.7).
\qed\medskip

\nin{\sl Proof  Theorem (3.9).} In outline the proof is the same as that of (2.4).  Fix $z\in E$. Hunt's switching identity \nref{\BG}{Thm.~VI(1.16)} implies that the duality between $X$ and $\hat X$ persists if each is killed at the first hitting time of $z$. Thus, writing $\EEE^f$ for $\int_E \EEE^x[\,\cdot\,] f(x)\,\pi(dx)$, etc., we have
$$
\EEE^f[g(X_t); t<T_z] =\hat \EEE^g[f(\hat X_t); t<\hat T_z],\qquad \forall f,g\in bp\EE.
\leqno(3.19)
$$
Integrating out $t>0$ and then setting $g\equiv1$ we obtain
$$
\EEE^f[T_z] =\hat \EEE^\pi\left[\int_0^{\hat T_z} f(\hat X_t)\,dt\right].
\leqno(3.20)
$$
Now integrate out $z$:
$$
\EEE^f[T_Z] =\hat \EEE^\pi\left[\int_0^{\hat T_Z} f(\hat X_t)\,dt\right],\qquad\forall f\in bp\EE.
\leqno(3.21)
$$
Note that the $\hat\PPP^\pi$-law of $\hat X_{\hat T_Z}$ is $\pi$. 

In view of Lemma (3.16)(b), if $K(x)=\infty$ for all $x$ there is nothing more to show.
So in the rest of the proof  we assume that $K(x)<\infty$ for all $x\in E$.

On the one hand, the right side of (3.21) is an invariant measure, by the argument used to prove Lemma (3.12). Because  $K$ is everywhere finite, the measure $K(x)\,\pi(dx)$ against which $f$ is being integrated to form the left side of (3.21)) is  $\sigma$-finite. Thus  Lemma (3.12) applies and we conclude that $\hat\kappa=\hat\EEE^\pi[\hat T_Z]<\infty$.
Moreover, still by   Lemma (3.12),   the right side of (3.21)  is equal to
$$
\pi(f)\hat\kappa.
\leqno(3.22)
$$
Varying $f$ in (3.21), and using (3.22), we see that
$$
\EEE^x[T_Z] = \hat\kappa,\qquad \pi\hbox{-a.e. } x\in E,
\leqno(3.23)
$$
proving (3.11). 

To prove (3.10) we are going to show that $K$ is finely lower semicontinuous. For each $z$, the terminal time $T_z$ is {\it exact}:
$$
\lim_{t\downarrow 0}T_z\comp\,\theta_t =T_z.
$$
Consequently, the bounded function
$$
\varphi^{\alpha}_z: x\mapsto \EEE^x[\exp(-\alpha T_z)]
\leqno(3.24)
$$
is  $\alpha$-excessive, hence finely continuous.  As noted by Meyer \reef{\M}, this implies that  the ``mixture'' 
$$
\varphi^\alpha: x\mapsto \EEE^x[\exp(-\alpha T_Z)],
\leqno(3.25)
$$
obtained by integrating out $z$, is likewise finely continuous. Consequently
$$
\EEE^x[T_Z] =\uparrow\lim_{\alpha\downarrow 0}{1-\varphi^\alpha(x)\over\alpha}
\leqno(3.26)
$$
is finely lower semicontinuous. The $\pi$-null  set $\{x: \EEE^x[T_Z]>\hat\kappa\}$ is therefore finely open, hence empty.  
It follows that $K(x)\le\hat\kappa$ for all $x\in E$.
\qed
\bigskip

\b{4. Eliminating The Exceptional Set}
\medskip

In this  section we refine the argument of section 3 to show that  $K(x)=\hat\kappa$ for all $x\in E$. Thus the exceptional set in  (3.11) is actually empty.

We begin with the
\medskip

\nin{\sl Proof of Proposition (3.7).}
Fix a state $y$. Suppose first that $y$ is {\it irregular$\,$}; that is, $\PPP^y[T_y>0]=1$. 
Because of (3.1), the process obtained by killing $X$ at time $T_y$ is transient, so there is a Borel function $q:E\to ]0,1]$ with
$$
\EEE^x\int_0^{T_y} q(X_t)\, dt\le 1,\qquad\forall x\in E.
\leqno(4.1)
$$
See \nref{\Gtrans}{Prop.~(2.2)}. We now apply Lemma (3.12) with $\mu$ the unit point mass at $y$ and $S=T_y$. The $\sigma$-finiteness condition is met because of the function $q$: $\EEE^y\int_0^{T_y}q(X_t)\,dt\le 1$.  It follows that $\EEE^y[T_y]<\infty$. 
Now let $x$ be another state, distinct from $y$.  Then by the strong Markov property,
$$
\eqalign{
 \infty> \EEE^y[T_y]
&\ge \EEE^y[T_y; T_x<T_y] =\EEE^y[T_x+T_y\comp\,\theta_{T_x}: T_x<T_y)]\cr
&\ge \EEE^y[T_y\comp\,\theta_{T_x}: T_x<T_y] = \PPP^y[T_x<T_y]\cdot\EEE^x[T_y].\cr   
}
\leqno(4.2)
$$
Clearly $\PPP^y[T_x<T_y]>0$ by the hypothesis $\PPP^y[T_x<\infty]=1$, because the first hit of $x$ must occur during an excursion away from $y$.
It follows that $\EEE^x[T_y]<\infty$ for 
 for all  $x\in E$ different from $y$, and we've  already shown that $\EEE^y[T_y]<\infty$.
 
The discussion when $y$ is a {\it regular} point ({\it i.e.}, $\PPP^y[T_y=0]=1$)  is similar, but in this case we use the 
 local time process $\{L^y_t: t\ge 0\}$ at $y$ and the associated It\^{o} excursion measure $n_y$. See \reef{\Gexc}, \reef{\FG} for background. 
The relevant fact is that 
$$
\EEE^y\left[\int_0^\infty e^{-\alpha t}\,d_tL^y_t\right]
={1\over \alpha c(y)+\alpha\int_0^\infty e^{-\alpha u}  n_y[T_y>u]\,du},
\leqno(4.3)
$$
where $c(y)\ge 0$ is a constant --- the stickiness of the process at $y$.  (Here the subscript $t$ in $d_tL_t^y$ indicates a differential with respect to $t$.)  
By the ergodic theorem  \nref{\F}{Thm.~6.1}, $\lim_{\alpha\to 0+}\alpha \EEE^y\left[\int_0^\infty e^{-\alpha t}\,d_tL^y_t\right]=\nu_{L^y}(1)$, the total mass of the Revuz measure of the continuous additive functional (CAF)  $L^y$, a positive finite constant. See \nref{\FG}{p.~421}. It follows from this and (4.3) that
$$
n_y[T_y]:=\int_\Omega T_y(\omega)\, n_y(d\omega) =\int_0^\infty n_y[T_y>u]\, du <\infty.
\leqno(4.4)
$$
Let  $x$ be another state, distinct from $y$. As before, the hypothesis $\PPP^y[T_x<\infty]=1$ implies that  $n_y[T_x<T_y]>0$. 
By the strong Markov property under $n_y$, \nref{\Gexc}{Thm.~(2.5)},
$$
\eqalign{
 \infty> n_y[T_y]
&\ge n_y[T_y; T_x<T_y] =n_y[T_x+T_y\comp\,\theta_{T_x}: T_x<T_y]\cr
&\ge n_y[T_y\comp\,\theta_{T_x}: T_x<T_y]= n_y[T_x<T_y]\cdot\EEE^x[T_y].\cr
}
\leqno(4.5)
$$
It follows that $\EEE^x[T_y]<\infty$. Trivially, $\EEE^y[T_y] =0<\infty$. \qed
\bigskip

In the remainder of this section we assume that $K$ is finite, the other case being covered completely by Theorem (3.9) because of Lemma (3.16)(b).  
In particular,  $\kappa = \hat\kappa<\infty$.

We now show that the expectation $\EEE^x[T_y]$ is not merely finite, but bounded as a function of $x$, for fixed $y$. 
First note that upon integrating out $y$ in (3.17) we obtain
$$
\EEE^x[T_z]\le K(x)+\EEE^\pi[T_z],\qquad\forall x,z\in E.
\leqno(4.6)
$$
Likewise, integrating out $x$ in (3.17):
$$
\EEE^\pi[T_z]\le\EEE^\pi[T_y]+\EEE^y[T_z]
\leqno(4.7)
$$

\proclaim{(4.8) Proposition} 
We have
$$
\EEE^\pi[T_z]<\infty,\qquad\forall z\in E.
\leqno(4.9)
$$
Consequently,  for each fixed $z\in E$,
$$
\EEE^x[T_z]\le C(z),\qquad\forall x\in E,
\leqno(4.10)
$$
where $C(z):=\hat\kappa + \EEE^\pi[T_z]<\infty$.
\endproclaim
\medskip

\nin{\sl Proof.}
Clearly
$$
\int_E\EEE^\pi[T_z]\,\pi(dz) =\kappa<\infty,
$$
and so 
$$
\EEE^\pi[T_z]<\infty,\qquad\pi\hbox{-a.e.} z\in E.
$$
From this,  (4.7), and Proposition (3.7), we deduce (4.9).
The estimate (4.10) now follows from (4.6) and (4.9) because $K(x)\le\hat\kappa$ for all $x$.
\qed
\medskip

The proof just given relies on the $\pi$-integrability of $z\mapsto\EEE^\pi[T_z]$.
Unlike $K$, however, this function need not be bounded.  See Example (6.1).

It can be shown that the uniform bound (4.10) is equivalent to $\lim_{\alpha\to 0+}\sup_{x\in E}|\alpha u^\alpha_{L^y}(x)-b(y)|=0$, for a certain constant $b(y)$, where $u^\alpha_{L^y}(x):=\EEE^x\int_0^\infty e^{-\alpha t} d_t L^y_t$ is the $\alpha$-potential of the local time $L^y$ used above for regular $y$ (and understood to be the simple counting process when $y$ is irregular). This kind of equivalence is familiar in the stability theory of Markov processes. See, for example, \nref{\MT}{Thm.~16.2.2} for discrete time and \nref{\DMT}{Sect.~6} for continuous time.

We need a lemma due to  Khas'minskii \nref{\K}{Lem.~3}; see also \nref{\CZ}{Lem.~3.7}. Recall that a  terminal time is a stopping time $T$ such that $T=t+T\comp\,\theta_t$ on $\{t<T\}$, for each $t>0$. Hitting times like $T_z$ are the prototypes of terminal times.

\proclaim{(4.11) Lemma} Let $T$ be a terminal time and suppose that $\EEE^x[T]\le C<\infty$ for all $x\in E$. Then
$$
\EEE^x[T^k]\le k!\cdot C^k,\qquad \forall x\in E,\quad k=2, 3, \ldots..
\leqno(4.12)
$$
\endproclaim

\proclaim{(4.13) Theorem} The Kemeny function $K$  is constant:  
$$
K(x) =\hat\kappa,\qquad\forall x\in E.
\leqno(4.14)
$$
\endproclaim

\nin{\sl Proof.} Let $C$ be as in (4.10) and assume for the moment that   
$$
\beta:=\int_E C(z)^2\,\pi(dz)<\infty. 
\leqno(4.15)
$$
By (4.10) and Lemma (4.11), we have
$$
\EEE^x[T_Z^2]\le 2\beta<\infty,\qquad\forall x\in E.
$$
Consequently, the elementary inequality
$$
\Big|{1-e^{-\alpha t}\over\alpha}-t\Big|\le {\alpha t^2\over 2},\qquad t\ge 0, \alpha>0,
$$
implies that the convergence of $[1-\varphi^\alpha(x)]/\alpha$ 
to $\EEE^x[T_Z]$ (used already in the proof of Theorem (3.9))  is uniform in $x$. Because $\varphi^\alpha(x) =\EEE^x[\exp(-\alpha T_Z)]$ is finely continuous, so is   $K$; thus  (3.11) implies (4.14).

To remove the extra hypothesis we proceed as follows. Choose a Borel function  $\psi:E\to ]0,1]$ such that $\pi(\psi\cdot C^2)<\infty$.  
Define $\psi_n(x):=\min(n\cdot\psi(x),1)$, and $\pi_n(dx):=[\pi(\psi_n)]^{-1}\psi_n(x)\,\pi(dx)$. Clearly $\psi_n$ increases pointwise to $1$, and $\pi_n$ converges setwise to $\pi$. 
Now use $\psi_n/\pi(\psi_n)$ to time change $X$ and $\hat X$, obtaining  Hunt processes $X^{(n)}$ and $\hat X^{(n)}$ in duality with respect to $\pi_n$. More precisely, with
$$
A^{(n)}_t:=[\pi(\psi_n)]^{-1}\int_0^t\psi_n(X_s)\,ds,\qquad t\ge 0,
$$
and
$$
\tau^{(n)}(t):=\inf\{s>0: A^{(n)}_s>t\},\qquad t\ge 0,
$$
we have
$$
X^{(n)}_t=X_{\tau^{(n)}(t)},\qquad t\ge 0,
$$
with the analogous definition for $\hat X^{(n)}$.   
Using the obvious notation, we have  $T^{(n)}_y =A^{(n)}_{T_y}$, and so
$$
\EEE^x[T^{(n)}_y]\le{C(y)\over\pi(\psi_n)},\qquad\forall x\in E,
$$
because $\psi_n(\cdot)\le 1$. Now this bound is square integrable with respect to $\pi_n(dy)$; that is, (4.15) holds for $X^{(n)}$. By the discussion of the preceding paragraph, $K_n(x)=\hat\kappa_n$ for all $x\in E$.  Evidently, $\lim_n\hat\kappa_n =\hat \kappa$ and $K(x) = \lim_n K_n(x)$ for all $x$, so $K$ is also a constant function.\qed
\bigskip

\vfill\eject

\b{5. All Points Regular}
\medskip

In this  section we  give another refinement of Theorem (3.9), based on a remarkable observation of   Eisenbaum and Kaspi \reef{\EK}.
From here on out we assume that  {\it all} points are regular:
$$
\PPP^x[T_x=0]=1,\qquad\forall x\in E.
\leqno(5.1)
$$
For example, $X$ could be a regular $1$-dimensional diffusion, or a continuous-time Markov chain on a countable state space $E$ with only stable states.

Because $X$ and $\hat X$ have the same semipolar sets \nref{\BG}{VI(1.19)}, the dual of (5.1) is automatically satisfied.

Because of (5.1),  $X$ admits a local time process $\{L^x_t: x\in E, t\ge 0\}$. For each $x$, $t\mapsto L^x_t$ is a CAF of $X$ that increases only on the visiting set $\{s\ge 0: X_s=x\}$. From \nref{\GK}{Thm.$\,$1}  we know  that this collection of CAFs can be chosen to be jointly measurable in $(t,x)$, and  normalized to serve as  occupation density with respect to $\pi$; that is, 
almost surely,
$$
\int_0^t f(X_s)\,ds=\int_E f(x)L^x_t\,\pi(dx),\qquad\forall t\ge 0, f\in bp\EE.
\leqno(5.2)
$$
We have
$$
u^\alpha(x,y)=\EEE^x\int_0^\infty e^{-\alpha t}\,d_tL_t^y,\qquad\forall x,y,\in E,
\leqno(5.3)
$$
along with the dual equality involving the dual local time process $\{\hat L^x_t: x\in E, t\ge 0\}$ for $\hat X$. 

Define, for $x,y,z\in E$,
$$
\eqalign{
    v_z(x,y)
    &:=\EEE^x[L^y_{T_z}]=\lim_{\alpha\to 0+}\EEE^x\int_0^{T_z} e^{-\alpha t}\,d_t L^y_t\cr
&=\lim_{\alpha\to 0+} \left\{u^\alpha(x,y) -\EEE^x[e^{-\alpha T_z}]\cdot u^\alpha (z,y)\right\}\cr
&=\lim_{\alpha\to 0+} \left\{u^\alpha(x,y) -u^\alpha (x,z)\cdot\hat\EEE^y[e^{-\alpha T_z}] \right\}\cr
&=\hat\EEE^y[\hat L^x_{\hat T_z}].\cr}
\leqno(5.4)
$$
(For the fourth equality above we have used Hunt's switching identity \nref{\BG}{Thm.~VI(1.16)}.) 
It should be noted that $v_z(x,y)<\infty$ for all $(x,y,z)\in E^3$. This follows from the fact that, by Proposition (3.7),
$$
\infty>\EEE^x[T_z] =\int_E v_z(x,y)\,\pi(dy),\qquad \forall x,z\in E,
$$
so that $y\mapsto v_z(x,y)$ is finite $\pi$-a.e. This function is excessive for $\hat X$ killed at time $T_z$ and so is everywhere finite on $E\setminus\{z\}$ because there are no non-empty polar sets for $\hat X$. And clearly $v_z(x,y) = 0 $ if $y=z$. 

By (5.2) and its dual, and Fubini,
$$
K(x) =\int_E\int_E v_z(x,y)\,\pi(dy)\,\pi(dz),\quad \hat K(y) =\int_E\int_E v_z(x,y)\,\pi(dx)\,\pi(dz).
$$
{\it Crucially,} 
$$
v_y(x,x) =v_x(y,y)=\EEE^x[T_y]+\EEE^y[T_x],\qquad\forall x,y\in E.
\leqno(5.5)
$$
The symmetry recorded as the first equality in (5.5) comes from \nref{\EK}{Lem.~2.1}.
It's clear from Proposition (3.7) that $\EEE^x[e^{-\alpha T_y}] = 1-\alpha\EEE^x[T_y] +o(\alpha)$ as $\alpha\to 0+$ and likewise for $\EEE^y[e^{-\alpha T_x}]$; sending $\alpha\to 0+$ in formula (7) in the proof of \nref{\EK}{Lem.~2.1} we find that $v_y(x,x)=\EEE^x[T_y]+\EEE^y[T_x]$, proving the second equality in (5.5). We have used here the fact that $\lim_{\alpha\to 0+}\alpha u^\alpha(x,x) =1$, as a consequence of the ergodic  theorem [13; Thm. 6.1] cited already just below (4.3), because the Revuz measure of $L^y$ is the unit point mass at $y$ by virtue of (5.2) and (5.3).

It will be useful later to know that, by \nref{\EK}{Lem.~2.2}, the formula 
$$
[d(x,y)]^2:=v_z(x,x)+v_z(y,y)-v_z(x,y)-v_z(y,x)\ge 0,\qquad x,y\in E,
\leqno(5.6)
$$
defines a metric $d$ on $E$, the remarkable thing being that the right side of (5.6) {\it doesn't depend on} $z$.  Furthermore, by  \nref{\EK}{Lem.~2.4}
$$
h(x,y):=[d(x,y)]^2=v_y(x,x),\qquad\forall x,y\in E.
\leqno(5.7)
$$
[The proof given in \reef{\EK} works for $x\not=y$ and $y\not=z$; the cases when  $x=z$ or $y=z$ are easy to verify directly.]

Define
$$
\gamma:=\int_E\int_E h(x,y)\,\pi(dx)\pi(dy).
\leqno(5.8)
$$
Of course, (5.5) implies that $\gamma =2\kappa =2\hat\kappa$, where $\kappa$ and $\hat\kappa$ are as defined in Theorem (3.9) and just below it.
Here is our improvement of Theorem (3.9), when all points are regular. 

\proclaim{(5.9) Theorem} 
 $K(x) =\hat K(x) =\gamma/2$ for all $x\in E$.
\endproclaim
\medskip

\nin{\sl Proof.} We begin by showing that
$$
v_z(x,y)+v_y(x,z) =h(y,z),\qquad\forall x,y,z\in E.\leqno(5.10)
$$
This is trivial if $y=z$ because all three terms vanish in that case. 
So fix $x$ and  $y\not=z$. By the strong Markov property,
$$
v_z(x,y)=\EEE^x[L^y_{T_z}]=\PPP^x[T_y<T_z]\cdot \EEE^y[L^y_{T_z}]=\PPP^x[T_y<T_z]\cdot v_z(y,y)
$$
and likewise for $v_y(x,z)$. 
Because $\PPP^x[T_y<T_z]+\PPP^x[T_z<T_y]=1$, the symmetry of $h$ 
implies
$$
\eqalign{
v_z(x,y)+v_y(x,z)
&=\PPP^x[T_y<T_z]v_z(y,y) + \PPP^x[T_z<T_y]v_y(z,z)\cr
&=\PPP^x[T_y<T_z]h(y,z) + \PPP^x[T_z<T_y]h(z,y)\cr
&=\left\{\PPP^x[T_y<T_z] + \PPP^x[T_z<T_y]\right\}h(y,z)\cr
&=h(y,z),\cr
}
$$
proving (5.10).
Now integrate (5.10) and  use the fact that the product measure $\pi\otimes\pi$ is invariant under $(y,z)\mapsto(z,y)$:
$$
\eqalign{
K(x)&=\int_E\int_E v_z(x,y)\pi(dy)\pi(dz)\cr
&={1\over 2}\left[\int_E\int_E v_z(x,y)\pi(dy)\pi(dz)+\int_E\int_E v_y(x,z)\pi(dy)\pi(dz)\right]\cr
&={1\over 2}\int_E\int_E h(y,z)\pi(dy)\pi(dz)
={\gamma\over 2},\cr
}
$$
establishing that $K$ is a constant function, with constant value equal to $\gamma/2$.
The dual equality $\hat K(x) =\gamma/2$ is proved the same way.
\qed
\bigskip

\b{(5.11) Remark} 
Given any probability distribution $\mu$ on $(E,\EE)$ with fine support $E$ (that is, $\mu$ charges each non-empty finely open subset of $E$), we can time change $X$ using the strictly increasing CAF $t\mapsto \int_E L^x_t\,\mu(dx)$ to obtain a second Hunt process $\tilde X$ with invariant distribution $\mu$ to which Theorem (5.9) applies. Writing $\tilde L^y_t$ for the associated local time process and $\tilde T_z$ for the hitting time of $z$ by $\tilde X$, one checks easily that $\tilde L^y_{\tilde T_z} =L^y_{T_z}$. Thus, by (5.4), the potential density $v_z(x,y)$ is {\it time-change invariant}; {i.e.}, $\tilde v_z(x,y) = v_z(x,y)$ for all $x,y,z\in E$.
The upshot is that
$$
\int_E\int_E v_z(x,y)\,\mu(dy)\mu(dz)
$$
is a constant function of $x\in E$, regardless of $\mu$. (This can also be proved directly by using (5.10) as in the proof of Theorem (5.9).) As $h$ is the square of a distance, we have
$$
h(x,y)\le 2h(x,x_0)+2 h(x_0,y),\qquad\forall x,y\in E,
$$
where $x_0$ is any fixed element of $E$.
Choose a Borel function  $\psi:E\to]0,1]$ such that 
$$
\int_E h(x,x_0)\psi(x)\,\pi(dx)<\infty.
$$
The above discussion applies to  $\mu(dx):=\psi(x)\pi(dx)$, and in this case  the Kemeny constant for $\tilde X$ is finite. That is, the original process $X$ is always time-change equivalent to  a Hunt process with finite Kemeny constant.
\bigskip

\b{(5.12) Remark}
Suppose we time change $X$ as in the preceding remark to obtain a second Hunt process $\tilde X$ with invariant distribution $\mu$. 
In view of the  discussion there,  (5.5)  shows that the {\it commute times} \nref{\LP}{p.~130}
$$
\EEE^x[T_y]+\EEE^y[T_x]
$$
are also time-change invariant. This is remarkable because  the separate mean hitting times clearly depend on the invariant distribution; for example,
$$
\EEE^x[T_z] =\int_E v_z(x,y)\,\pi(dy).
$$
\medskip

\b{(5.13) Remark}  
 The formula $\kappa =\gamma/2$  is consistent with the formula reported in \nref{\Pin}{Rem.$\,$3}, because $h(x,y) =2|S(x)-S(y)|$ when $X$ is a diffusion on $\RRR$ with scale $S$ as in \reef{\Pin}. Pinsky shows that $\gamma<\infty$ if and only if the boundary points $-\infty$ and $+\infty$ are both  entrance boundaries. One might conjecture that, in the context of sections 3--5 of this paper,  the appropriate boundary to consider is $\overline E\setminus E$, where $\overline E$ is  a  Ray-Knight compactification of $E$. Then, in this the best of all possible worlds,  the condition  $\gamma<\infty$ would be equivalent to the ``Ray space'' of $X$ coinciding with $\overline E$. See \nref{\GS}{Def.$\,$4.5}. We hope to return to this question in a future paper, but for now see Example (6.7).
\bigskip

\b{(5.14) Remark} Suppose now that $X$ is reversible  with respect to $\pi$ ({\it i.e.},  $X=\hat X$).
 Fix distinct $y$ and $z$. Then $u(x):=\PPP^x[T_y<T_z]$ is the condenser potential for the pair $\{y\},\{z\}$; see \reef{\CG} for a probabilistic discussion of these matters. The function $u$ is the  potential, with respect to $X$ killed at time $T_z$, of a point mass at $y$. That is,
$$
\PPP^x[T_y<T_z]= c(y,z)\cdot v_z(x,y),
$$
for some constant $c(y,z)\in]0,\infty[$. Using the  strong Markov property at time $T_y$ as before,
$$
v_z(x,y)=\PPP^x[T_y<T_z]\cdot \EEE^y[L^y_{T_z}]=c(y,z)\cdot u(x)\cdot h(y,z).
$$
It follows that
$$
h(y,z) = 1/c(y,z).
\leqno(5.15)
$$
Let $(\EE,\FF_e)$ denote the (extended) Dirichlet form associated with $X$. It is known that $u$ is the solution of the optimization problem
$$
\inf\{ \EE(w,w): w\in\FF_e, w(y) =1, w(z) =0\};
$$
see \nref{\CFS}{Sect.$\,$3}. As such, $c(y,z)$ is the Dirichlet energy $\EE(u,u)$ of $u$, commonly referred to as the {\it effective resistance} between $y$ and $z$. And then $h$ is the effective resistance distance between $y$ and $z$. This identification of $h$ may prove useful in computing $\gamma$ (and thereby $\kappa$) in specific situations.

\bigskip

\b{(5.16) Remark} An interesting consequence of Theorem (5.9) is the formula  
$$
\EEE^\pi[L^x_{T_Z}] =\int_E\int_E v_z(y,x)\,\pi(dy)\pi(dz) = \hat \EEE^x[\hat T_{Z}]=\gamma/2 = \kappa,\qquad\forall x\in E. 
\leqno(5.17)
$$
\bigskip

\b{6. Examples}
\medskip

\b{(6.1) Example}  Let $(B_t)_{t\ge 0}$ be standard Brownian motion in the punctured ball $\{x\in\RRR^3: 0<|x|\le 1\}$, reflected at the bounding sphere.
The radial part $X_t:=|B_t|$ is then a $3$-dimensional Bessel process with state space $E=]0,1]$, with reflection at $1$.
 The infinitesimal  generator of $X$ is
$$
\LL f(x) = {1\over 2}f''(x) +{1\over x}f'(x),\qquad 0<x\le 1,
\leqno(6.2)
$$
with boundary condition $f'(1) =0$. This process is positive recurrent, with stationary distribution
$$
\pi(dx) = 3x^2\,dx,
$$
(coincident with the speed measure of $X$) and associated scale function
$$
S(x) = -{2\over 3x},\qquad 0<x\le 1.
$$
 Fix $z\in ]0,1]$. Clearly $v_z(x,y)$ vanishes if $x$ and $y$ are on opposite sides of $z$. For $x,y\in[z,1]$ we have
 $$
 v_z(x,y) =[S(x)-S(z)]\wedge[S(y)-S(z)],
 $$
 resulting in
 $$
 \EEE^x[T_z] =\int_z^1 v_z(x,y)\,\pi(dy) ={z^2-x^2\over 3}+{2\over 3}\cdot\left({1\over z}-{1\over x}\right),\qquad 0<z\le x\le 1.
 \leqno(6.3)
 $$
 Meanwhile, for $x,y\in ]0,z]$  we have
 $$
  v_z(x,y) =[S(z)-S(x)]\wedge[S(z)-S(y)],
  $$
  resulting in 
  $$
  \EEE^x[T_z] = \int_0^z v_z(x,y)\,\pi(y) = {z^2-x^2\over 3},\qquad 0<x\le z.
  \leqno(6.4)
  $$
Integrating: 
$$
\EEE^x[T_Z] ={1\over 5},\qquad\forall x\in ]0,1],
$$
constant as expected. Notice that $h(x,y) = v_y(x,x) = |S(x) -S(y)|$, and therefore
$$
\gamma
=\int_0^1\int_0^1 |S(u)-S(v)|\,\pi(du)\,\pi(dv)
={2\over 5},
$$
consistent with Theorem (5.9). Incidentally, formulas (6.3) and (6.4) are consistent with (5.5).

One checks that
$$
\EEE^\pi[T_z] = {2\over 3z}-{6\over 5}+{2z^2\over 3},\qquad 0<z\le 1,
$$
which grows without bound as $z\downarrow 0$. This confirms the remark made just after the proof of Proposition (4.8), 
\bigskip

\b{(6.5) Example} Now let $X$ be the real-valued Ornstein-Uhlenbeck process, with generator
$$
\LL f(x) = {1\over 2}f''(x) -{x\over 2} f'(x),\qquad x\in\RRR.
\leqno(6.6)
$$
The process $X$ is positive recurrent with invariant distribution $\pi$ (normalized speed measure) equal to  the standard normal distribution. The associated scale function, chosen so that $S(0)=0$, is given by
$$
S(x) =\sqrt{8\pi}\int_0^x e^{u^2/2}\,du,\qquad x\in \RRR.
$$
As before, $h(x,y) =|S(x) -S(y)|$, so the rapid growth of $S$ implies that $\gamma=\infty$, so $K(x) =\infty$ for all $x\in\RRR$, by Theorem (5.9). Of course this follows from \nref{\Pin}{Thm.~1.1} because the boundaries $\pm\infty$ for the O.U. process are not entrance boundaries.

\bigskip

\b{(6.7) Example}  We now take $X$ to be the skew Brownian motion of It\^o-McKean \nref{\IM}{pp.$\,$115--116}; see also \reef{\HS}.  
This is a (reversible) diffusion on $[-1,1]$  described intuitively as follows. Fix $0<\alpha<1$. Let $Y$ be (standard) Brownian motion on $[0,1]$ with reflecting barriers at $0$ and $1$.  Change the sign of each excursion of $Y$ away from $0$ with probability $1-\alpha$; the various sign changes are understood to be independent of each other and of $Y$. The resulting process is a regular diffusion on $[-1,1]$ with reflecting barriers at $\pm 1$, scale function
$$
S(x) =
\cases{
{2x/\alpha},& $0\le x\le 1$,\cr 
{2x/ (1-\alpha)},& $-1\le x\le 0$,\cr
}
\leqno(6.8)
$$
and speed measure (and invariant distribution) 
$$
\pi(dx)/dx =
\cases{
\alpha, & $0< x\le 1$,\cr
1-\alpha, & $-1\le x< 0$.\cr
}
\leqno(6.9)
$$
We have normalized $\pi$ to be a probability distribution. Given  the choice $S(0)=0$, this forces the formula for $S$.
As in Example (6.1) we use $S$ and $\pi$ to compute 
$$
\kappa =4/3
\leqno(6.10)
$$
for all $-1<\alpha<1$. Notice that $\kappa =1/3$ in the special cases $\alpha =\pm1$ (reflected Brownian motion on $[0,1]$ or $[-1,0]$). 
This discontinuity in $\kappa$  at $\alpha=\pm 1$ might be expected, but I have no simple explanation for why $\kappa$ doesn't otherwise depend on $\alpha$.

More generally, for $n\ge 2$, let $X$ be the Walsh Brownian motion \reef{\W} on the rimless wheel with $n$ spokes, in $\RRR^2$:
$$
E=\{0\}\cup\left[\cup_{j=1}^n\{re^{i\theta_j} : 0<r\le 1\}\right],
$$
where $2\pi > \theta_1>\theta_2>\cdots >\theta_n>0$, with reflection at the outer endpoints $e^{i\theta_j}$, $j=1,2,\ldots,n$. See \reef{\FK} for a discussion of this process and further references. On each spoke, the process $X$ moves like Brownian motion reflected at the outer endpoint of the spoke.  The excursions of $X$ from the origin are performed along the various spokes (labelled by the angles $\theta_1,\ldots,\theta_n$); the spoke for a given excursion is chosen at random using probabilities $p_1,p_2,\ldots,p_n$ (with $p_j>0$ for each $j$ and $\sum_{j=1}^n p_j=1$), the choices for different excursions being mutually independent. Owing to the tree-like nature of $E$, one can show that $h(x,y) =|S(x)-S(y)|$, where $S$ is a sort of scale function, given on the $j^{\rm th}$ spoke by
$$
S(re^{i\theta_j} )= 2r/p_j,\qquad 0< r\le 1.
$$
Moreover, on spoke $j$, $\pi$ is just $p_j$ times Lebesgue measure. 
A straightforward calculation reveals that
$$
\kappa=n-{2\over 3}. 
\leqno(6.11)
$$
Like skew Brownian motion (the case $n=2$),  $\kappa$ doesn't depend on the distribution $p=\{p_j\}_{j=1}^n$.  We can even allow  a countable infinity of spokes $\{\theta_j\}_{j=1}^\infty$, with spoke weights given by some  $p_j>0$ with $\sum_{j=1}^\infty p_j=1$. If we   time change to sample this motion only when it is visiting the spokes with angles in $\{\theta_1,\ldots,\theta_n\}$, the result is  a Walsh Brownian motion with $n$ spokes.  It follows from this and  (6.11) that Kemeny's constant for the infinite-spoke motion is infinite.

This example lends support to a conjecture voiced in Remark (5.13). Let us be definite and choose $\theta_n=1/n$, $n=1,2,\ldots$. The relative Euclidean topology on $E$ coincides with the  topology on $E$ determined by the evident ``tree distance''. The Euclidean closure of $E$ is
$$
\overline E=E\cup \left(]0,1]\times\{0\}\right).
$$
This can serve as the state space for a Ray-Knight compactification of $X$. When started in a point $(r,0)$ of $\overline E\setminus E$, the behavior of the associated Ray process $\overline X$ can be described as follows: $\overline X_t=(B_t,0)$, where $B$ is a Brownian motion on $]0,1]$ started at $r\in ]0,1]$, with reflection at $1$, until it hits the origin; after that $\overline X$  moves like an  independent copy of $X$ started at the origin. Thus the boundary $\overline E\setminus E$ is disjoint from the Ray entrance space for $X$ (namely, $E$ itself). 
\bigskip

\b{(6.12) Example} This final example is a pure-jump process. Fix $\alpha\in]1,2]$, and let $Y=(Y_t)_{t\ge 0}$ be a symmetric $\alpha$-stable L\'evy process with $Y_0=0$ and 
$$
\EEE[e^{i\lambda Y_t}] =e^{-tC|\lambda|^\alpha},\qquad t\ge 0, \lambda\in\RRR,
$$
where $C>0$. Let $E$ denote the unit circle in the complex plane, and  wrap $Y$ around  $E$ to obtain the $E$-valued process
$$
X_t:=\exp(iu+iY_t),\qquad t\ge 0.
$$
Here $X_0=x:=e^{iu}$, where $0\le u<2\pi$. 
This is a positive recurrent, self-dual Hunt process with state space $E$ and invariant distribution the uniform law on $E$. 

For the moment   take $1<\alpha<2$.  If $y\in E$, say $y=e^{iv}$ with $0\le v<2\pi$, then
$$
T_y(X) = T_{v-u+2\pi\ZZZ}(Y),
$$
where $T_y(X)$ denotes the hitting time of $\{y\}$ by $X$ and $T_A(Y)$ the hitting time of $A$ by $Y$.
We know from \nref{\Iso}{Thm.~1(b)} that  
$$
\EEE[T_{v-u+2\pi\ZZZ}(Y) ]=\sum_{n\in\ZZZ,n\not=0}{1\over C|n|^\alpha}(1-e^{ -in (v-u)})={2\over C}\sum_{n=1}^\infty{1-\cos(n(v-u))\over n^\alpha}.
$$
Thus
$$
\EEE^x[T_y(X)] 
={2\over C}\sum_{n=1}^\infty {1-\cos(n(v-u))\over n^\alpha}.
$$
Because the integral of the cosine over a complete  period vanishes, we find that
$$
\eqalign{
K(x) 
&= \int_E \EEE^x[T_y(X)]\,\pi(dy) = {1\over 2\pi}\int_0^{2\pi} {2\over C}\sum_{n=1}^\infty {1-\cos(n(u-v))\over n^\alpha}\,dv\cr
&={2\over C}\sum_{n=1}^\infty {1 \over n^\alpha}<\infty.\cr
}
\leqno(6.13)
$$
As might be expected, because the symmetric Cauchy process ($\alpha =1$) doesn't hit points, this expression blows up when  $\alpha\downarrow 1$. 

On the other hand, if we take $\alpha =2$ and $C=1/2$, then $Y$ is  a standard Brownian motion. The hitting time of $y$ by $X$ (with $X_0=x$) has the same law as the exit time from an interval with endpoints $ |v-u|$ and $-(2\pi -|v-u|)$ by $Y$ (remember $Y_0=0$). As this exit time is well known to have mean $|v-u|\cdot(2\pi-|v-u|)$, the Kemeny constant for the wrapped Brownian motion is easily seen to be $2\pi^2/3$.
 This coincides with the limit as $\alpha\uparrow 2$ in (6.13), as expected.


\bigskip\bigskip

\baselineskip=13pt
\frenchspacing

\centerline{\bf References}
\medskip

\itemitem{[\AF]} Aldous, D. and Fill, J.:  {\it Reversible Markov chains and random walks on graphs,} (2002).\hfill\break
{\tt http://www.stat.berkeley.edu/users/aldous/RWG/book.html}
\smallskip

\itemitem{[\AH]}
Angel, O. and Holmes, M.: Kemeny's constant for infinite DTMCs is infinite (Letter to the editor), {\it J. Appl. Prob.} {\bf  56} (2019) 1269--1270.
\smallskip

\itemitem{[\AKDR]}
Az\'ema, J., Kaplan-Duflo, M., Revuz, D.: Mesure invariante sur les classes r\'ecurrentes des processus de Markov, {\it Z.   Warsch. Verw. Geb.} {\bf 8} (1967) 157--181.
\smallskip

\itemitem{[\BHLMT]}
Bini, D., Hunter, J.J., Latouche, G., Meini, B., Taylor, P.:
Why is Kemeny's constant a constant? {\it J. Appl. Probab.} {\bf 55} (2018) 1025--1036. 
\smallskip

\itemitem{[\BG]}
Blumenthal, R.M. and Getoor, R.K.: {\it Markov Processes and Potential Theory.} Academic Press, New York, 1968.
\smallskip

\itemitem{[\CFS]}
Chen, Z.-Q., Fitzsimmons, P.~J.,  and Song, R.: { Crossing estimates for
  symmetric {M}arkov processes}, {\it Probab. Theory Related Fields}  {\bf 120}
  (2001) 68--84. 
\smallskip

\itemitem{[\C]}
Chung, K.L.: Contributions to the theory of Markov chains, II, {\it Trans. Amer. Math. Soc.} {\bf 76} (1954) 397--419. 
\smallskip

\itemitem{[\CG]}
Chung, K.L. and Getoor, R.K.: The condenser problem, {\it Ann. Probab.} {\bf 5} (1977) 82--86.
\smallskip

\itemitem{[\CZ]}
Chung, K.L. and Zhao, Z.:  {\it From Brownian Motion to Schr\" odinger's Equation}, Springer, Berlin, 1995.
\smallskip

\itemitem{[\Dev]}
Devriendt, K.: Kemeny's constant and the Lemoine point of a simplex, {\it Elec. J. Lin. Alg.} {\bf 40} (2024) 766--773. 
\smallskip

\itemitem{[\DMT]}
Down, D., Meyn, S. P., and Tweedie, R. L.: {Exponential and uniform ergodicity of {M}arkov processes}, {\it  Ann. Probab.} {\bf 23} (1995) 1671--1691. 
\smallskip

\itemitem{[\D]}
Doyle, P.: The Kemeny constant of a Markov chain, 10 pages, 2009, 
{\tt https://arxiv.org/pdf/0909.2636}
\smallskip

\itemitem{[\EK]}
Eisenbaum, N. and Kaspi, H.: On the continuity of local times of Borel right Markov processes, {\it  Ann.  Probab.} {\bf 35} (2007)  915--934.
\smallskip

\itemitem{[\F]}
Fitzsimmons, P.J.: The quasi-sure ratio ergodic theorem, {\it Ann. Inst. H. PoincarŽ Probab. Statist.} {\bf 34} (1998) 385--405. 
\smallskip

\itemitem{[\FG]}
Fitzsimmons, P.J. and Getoor, R.K.:
Excursion theory revisited, {\it Illinois J. Math.} {\bf  50}  (2006) 413--437. 
\smallskip

\itemitem{[\FK]}
Fitzsimmons, P.J. and Kuter, K.E.: Harmonic functions on Walsh's Brownian motion, {\it Stochastic Processes and their Applications}  {\bf  124} (2014) 2228--2248.
\smallskip

\itemitem{[\Gexc]}
Getoor, R.K.: Excursions of a Markov process, {\it Ann. Probab.} {\bf  7}  (1979) 244--266. 
\smallskip

\itemitem{[\Gtrans]}
Getoor, R.K.: Transience and recurrence of Markov processes,  {\it S\'eminaire de Probabilit\'es, XIV}, pp. 397--409,
Lecture Notes in Math. {\bf 784}, Springer, Berlin, 1980. 
\smallskip

\itemitem{[\Gexcess]}
Getoor, R.K.: {\it Excessive Measures,} Birkh\" auser, Boston, 1990.
\smallskip

\itemitem{[\GK]}Getoor, R.K. and Kesten H.: Continuity of local times for Markov processes, {\it Compositio Mathematica} {\bf 24} (1972) 277--303.
\smallskip

\itemitem{[\GS]}
Getoor, R.K. and Sharpe, M.J.: The Ray space of a right process, 
{\it Annales de l'Institut Fourier} {\bf  25} (1975)   207--233.
\smallskip

\itemitem{[\H]}
Harris, T.E.: The existence of stationary measures for certain Markov processes, {\it  Proc. Third Berkeley Symposium on Mathematical Statistics and Probability, 1954--1955}, vol. II, pp. 113--124, UC Press, Berkeley, 1956. 
\smallskip

\itemitem{[\HS]}
Harrison, J.M. and Shepp, L.A.: On skew Brownian motion, {\it Ann. Probab.} {\bf 9} (1981) 309--313.
\smallskip
\itemitem{[\HY]}
He, P. and Ying, J.: Fine irreducibility and uniqueness of stationary distribution, {\it Osaka J. Math.} {\bf 50} (2013)  417--423.
\smallskip

\itemitem{[\HunMix]}
Hunter, J.J.: Mixing times with applications to perturbed Markov chains,  {\it Linear Algebra Appl.} {\bf  417} (2006)  108--123. 
\smallskip

\itemitem{[\Hun]}
Hunter, J.J.: The role of Kemeny's constant in properties of Markov chains, {\it Comm. Statist. Theory Methods} {\bf  43} (2014)  1309--1321. 
\smallskip

\itemitem{[\KM]}
Kaspi, H. and Mandelbaum, A.: On Harris recurrence in continuous time, {\it Math. Oper. Res.} {\bf 19} (1994) 211--222.
\smallskip

\itemitem{[\KS]} 
Kemeny, J.G. and J.L. Snell: {\it Finite Markov Chains.} Van Nostrand,  Princeton-Toronto-London-New York, 1960.
\smallskip

\itemitem{[\K]}
Khas'minskii, R. Z.: On positive solutions of the equation ${\AA}u+Vu=0$, {\it Theor. Probability Appl.} {\bf  4} (1959) 309--318. 
\smallskip

\itemitem{[\LeLo]}
Levene, M. and  Loizou, G.: Kemeny's constant and the random surfer, {\it Am. Math. Monthly} {\bf 109}  741--745 (2002). 
\smallskip

\itemitem{[\LiLy]}
Liu, Y. and Lyu, F.: Kemeny's constant for countable Markov chains, {\it Linear Alg. Appl.} {\bf 604} (2020) 425--440.
\smallskip

\itemitem{[\M]}
Meyer, P.-A.: Une remarque sur la topologie fine, {\it S\'eminaire de Probabilit\'es, XIX}, p.$\,$176,
Lecture Notes in Math. {\bf 1123}, Springer, Berlin, 1985. 
\smallskip

\itemitem{[\MT]}
Meyn, S. and  Tweedie, R.L.: {\it Markov Chains and Stochastic Stability} (2nd edition), Cambridge, 2009.
\smallskip

\itemitem{[\Pin]}
Pinsky, R.: Kemeny's constant for one-dimensional diffusions, {\it Electron. Commun. Probab.} {\bf 24} (2019), no. 36, 5 pages.
\smallskip

\itemitem{[\Sh]}
Sharpe. M.J.: {\it General Theory of Markov Processes},  Academic Press,  Boston, 1988. 
\smallskip

\itemitem{[\Si]}
Silverstein, M.L.: Classification of coharmonic and coinvariant functions for a L\'evy process, {\it Ann. Probab.} {\bf 8} (1980) 539--575. 
\smallskip

\itemitem{[\W]}
Walsh, J.B.:  A diffusion with a discontinuous local time, {\it Ast\'erisque} {\bf 52--53} (1978) 37--45.
\smallskip

\end